\documentclass[12pt,english]{article}
\usepackage[latin1]{inputenc}
\usepackage{amsmath}
\usepackage{graphicx}
\usepackage{amssymb}

\makeatletter
 \usepackage{amsthm}
 \theoremstyle{plain}    
 \newtheorem{statement}{Statement} 
 \theoremstyle{plain}
 \theoremstyle{definition}
 \newtheorem{definition}[statement]{Definition}
 \theoremstyle{plain}    
 \newtheorem{theorem}[statement]{Theorem} 
 \theoremstyle{plain}    
 \newtheorem{lemma}[statement]{Lemma}
 \theoremstyle{plain}    
 \newtheorem{conjecture}[statement]{Conjecture}

\usepackage{color}
\usepackage[letterpaper,text={6.5in,9in}]{geometry}

\usepackage{babel}
\makeatother
\begin{document}

\title{Fire Containment in Grids of Dimension Three and Higher}

\author{Mike Develin\\
 American Institute of Mathematics\\
 360 Portage Avenue\\
 Palo Alto, CA 94306\\
 \texttt{develin@post.harvard.edu}\\
 \and Stephen G. Hartke%
\footnote{The second author has been supported by a National Science Foundation
Graduate Research Fellowship, a National Defense Science and Engineering
Graduate Fellowship, a DIMACS Graduate Research Award, and National
Science Foundation grants EIA-0205116, DBI-9982983, and SBR-9709134
through DIMACS. Part of this work was done while the two authors were
visiting the Research Experience for Undergraduates at the University
of Minnesota-Duluth, partially supported by National Science Foundation
grant DMS-0337448.%
}\\
 Department of Mathematics\\
 University of Illinois\\
 Urbana, IL 61801\\
 \texttt{hartke@math.uiuc.edu}}

\date{September 9, 2004}

\maketitle
\begin{abstract}
We consider a deterministic discrete-time model of fire spread introduced
by Hartnell~{[}1995{]} and the problem of minimizing the number of
burnt vertices when deploying a limited number of firefighters per
timestep. While only two firefighters per timestep are needed in the
two dimensional lattice to contain any outbreak, we prove a conjecture
of Wang and Moeller~{[}2002{]} that $2d-1$ firefighters per timestep
are needed to contain a fire outbreak starting at a single vertex
in the $d$-dimensional square lattice for $d\geq3$; we also prove
that in the $d$-dimensional lattice, $d\geq3$, for each positive
integer $f$ there is some outbreak of fire such that $f$ firefighters
per timestep are insufficient to contain the outbreak. We prove another
conjecture of Wang and Moeller that the proportion of elements in
the three-dimensional grid $P_{n}\times P_{n}\times P_{n}$ which
can be saved with one firefighter per time step when an outbreak starts
at one vertex goes to $0$ as $n$ gets large. Finally, we use integer
programming to prove results about the minimum number of timesteps
needed and minimum number of burnt vertices when containing a fire
outbreak in the two dimensional square lattice with two firefighters
per timestep.

\medskip{}
\end{abstract}

\section{Introduction}

Hartnell~\cite{Hartnell} introduced a deterministic discrete-time
model of fire spread on a graph $G$ and considered how firefighters
can act to stop a fire outbreak. In this model, an outbreak of fire
starts at a set of root vertices of $G$ at time $t=0$. In response,
firefighters are placed at the vertices $a_{1,1},a_{1,2},\ldots,a_{1,c_{1}}$
at time $t=1$, where the firefighters defend or protect each vertex
from the spreading fire. The fire then spreads from burning vertices
to non-defended neighbors. Firefighters are again deployed to defend
the vertices $a_{2,1},\ldots,a_{2,c_{2}}$ at time $t=2$ (the vertices
$a_{1,1},a_{1,2},\ldots,a_{1,c_{1}}$ remain defended), and the fire
spreads again. The process continues until the fire can no longer
spread. We say that the fire outbreak is \emph{contained} after $t$
time steps if there is some finite time $t$ such that after the disease
spreads during time $t$, only a finite number of vertices are burnt
and the disease can no longer spread. The motivating question is to
find an optimal sequence of defended vertices that minimizes the total
number of burnt vertices.

The fire spread model is also relevant in epidemiology. Traditionally,
epidemiological models assume that the population being studied is
well-mixed in the sense that any pair of individuals are just as likely
to come in contact and transmit a disease as any other. However, recently
epidemiologists have attempted to incorporate spatial information
into their models \cite{ReadKeeling,BurkeEpstein,EpiSIMS}. The model
of fire spread presented above can be considered as modeling a perfectly
contagious disease with no cure, where vertices adjacent to infected
vertices become infected at every discrete time step and, once infected,
remain infected from then on. The response allowed is only a limited
number of vaccinations of non-infected vertices per time step. The
limited number of vaccinations is particularly relevant to real-world
situations because of limited availability of the vaccine or limited
numbers of health personnel to administer the vaccine. The main question
is of course to minimize the total number of infected vertices.

In this work we study fire containment on square grids. Grids are
a natural class of graphs on which to consider both fire and disease
spread since they are often used to represent geographic areas. Both
Wang and Moeller~\cite{WangMoeller} and MacGillivray and Wang~\cite{MacGillivrayWang}
studied grids to find algorithms for containment. Wang and Moeller
showed that two firefighters per time step is sufficient to contain
a fire outbreak in a two dimensional square grid, and conjectured
that $2d-1$ firefighters are necessary to contain a fire outbreak
in a $d$ dimensional square grid. We prove this conjecture in section~\ref{sec:ThreePlusDgrids}.
Fogarty~\cite{Fogarty} showed that two firefighters suffice in the
two dimensional square lattice to contain any finite outbreak of fire
where an arbitrarily large but finite number of vertices are initially
on fire. However, we prove that for any fixed number $f$ of firefighters,
there is a finite outbreak of fire in which $f$ firefighters per
time step are insufficient to contain the outbreak. We also prove
the conjecture of Wang and Moeller that the proportion of elements
in the three-dimensional grid $P_{n}\times P_{n}\times P_{n}$ which
can be saved by using one firefighter per time step when an outbreak
at one vertex occurs goes to $0$ as $n$ gets large.

In section~\ref{sec:TwoDgrid} we provide an alternate proof using
integer programming of Wang and Moeller's result that the minimum
number of time steps needed to contain a fire outbreak in a two dimensional
square grid when using two firefighters per time step is $8$. We
also use this technique to prove that the minimum number of burned
vertices in such an outbreak is $18$.

We use the following terminology to describe the fire spread and firefighter
response. During the $t^{\textrm{th}}$ time step for $t>0$, firefighters
are deployed and then the fire spreads. If we describe the state of
vertices at the \emph{beginning of the $t^{\textrm{th}}$ time step},
we mean \emph{before} the firefighters are deployed during the $t^{\textrm{th}}$
time step. If we describe the state of vertices at the \emph{end}
\emph{of the $t^{\textrm{th}}$ time step}, or equivalently, at the
\emph{end of $t$} \emph{time steps,} we mean \emph{after} the fire
has spread during the $t^{\textrm{th}}$ time step. A firefighter
may defend neither a burnt vertex nor a previously defended vertex.
Once fire has spread to a vertex $v$, we say that $v$ is a \emph{burnt}
vertex. After being burnt or defended, a vertex remains in that state
until the process ends. In addition to the burnt and defended vertices,
we say that a vertex $v$ is \emph{saved} at the end of the $t^{\textrm{th}}$
time step if there is no path from $v$ to the root consisting only
of burnt and non-defended vertices at the end of the $t^{\textrm{th}}$
time step.

We consider the infinite $d$-dimensional square grids $\mathbb{L}^{d}$.
The vertices of $\mathbb{L}^{d}$ are the points of $\mathbb{R}^{d}$
with integer coordinates, and $x$ is adjacent to $y$ if and only
if $x$ is distance $1$ from $y$ in the usual Euclidean $\ell_{2}$
metric.

\section{\label{sec:ThreePlusDgrids}Three and Higher Dimensional Square Grids}

Wang and Moeller proved in \cite{WangMoeller} that an outbreak starting
at a single point in a regular graph of degree $r$ can be contained
with if $r-1$ firefighters can be deployed per time step. Specifically,
for the $d$ dimensional square grid $\mathbb{L}^{d}$, $2d-1$ firefighters
suffice to contain an outbreak starting at a single point. They conjectured
that this bound is tight, and we present a proof of this conjecture
here.

Wang and Moeller observed that at least two firefighters per time
step are needed for containment in $\mathbb{L}^{3}$, and Fogarty
showed in \cite{Fogarty} that at least three firefighters per time
step are needed to contain the outbreak. Her main theorem involves
a {}``Hall-type condition'' for the graph, which provides a lower
bound for how fast the fire can spread. The theorem considers the
\emph{front} of the fire, which is the set of burnt vertices farthest
from the root. The theorem states that if this front grows quickly
(\emph{i.e.}, by at least $f$) regardless of its precise shape, then
it cannot be contained by deploying $f$ firefighters per time step.
Theorem~\ref{thm:HallTypeCondition} strengthens Fogarty's theorem
by considering initial growth of the fire that is faster than $f$
so that the fire reaches a {}``critical mass\char`\"{} and can sustain
growth of the front by at least $f$ from that point onward.

First we state some definitions.

\begin{definition}
\label{def:DkAndrtAndBt}Let $D_{k}$ denote the set of vertices in
a rooted graph $G$ that are distance $k$ from the root vertex $r$.
Let $r_{k}$ denote the number of firefighters in $D_{k+1},D_{k+2},\ldots$
at the end of the $k^{\textrm{th}}$ time step. These firefighters
can be thought of as {}``reserve'' firefighters since they are not
adjacent to the fire when deployed. We define $r_{0}$ to be $0$.
Let $B_{k}\subseteq D_{k}$ denote the number of burned vertices in
$D_{k}$ at the end of the $k^{\textrm{th}}$ time step. 
\end{definition}
\begin{theorem}
\label{thm:HallTypeCondition} Let $G$ be a rooted graph, $h$ a
positive integer, and $a_{0},a_{1},\ldots,a_{h}$ positive integers
each at least $f$ such that the following holds: 
\begin{enumerate}
\item \label{enu:ConditionKequals0} Every $A\subseteq D_{0}$, $A\neq\emptyset$,
satisfies $\left|N(A)\cap D_{1}\right|\geq\left|A\right|+a_{0}$. 
\item \label{enu:Condition1KH} For $1\leq k\leq h$, every $A\subseteq D_{k}$
where $\left|A\right|\geq1+\sum_{i=0}^{k-1}(a_{i}-f)$ satisfies $\left|N(A)\cap D_{k+1}\right|\geq\left|A\right|+a_{k}$. 
\item \label{enu:ConditionKgtH} For $k>h$, every $A\subseteq D_{k}$ such
that $\left|A\right|\geq1+\sum_{i=0}^{h}(a_{i}-f)$ satisfies $\left|N(A)\cap D_{k+1}\right|\geq\left|A\right|+f$. 
\end{enumerate}
Suppose that at most $f$ firefighters per time step are deployed.
Then\begin{equation}
\left|B_{n}\right|\geq\begin{cases}
1 & \textrm{if $n=0$,}\\
1+r_{n}+\sum_{i=0}^{n-1}(a_{i}-f) & \textrm{if $1\leq n\leq h+1$,}\\
1+r_{n}+\sum_{i=0}^{h}(a_{i}-f) & \textrm{if $n>h+1$,}\end{cases}\label{eq:HallThmConclusion}\end{equation}
 regardless of the sequence of firefighter placements. Specifically,
$f$ firefighters per time step are insufficient to contain an outbreak
that starts at the root vertex. 
\end{theorem}
\begin{proof}
Let $p_{n+1}$ denote the number of firefighters placed in $D_{n+1}$
at time $n+1$, and let $p_{\leq n}$ denote the number of reserve
firefighters placed in $D_{n+1}$ during time steps $1,\ldots,n$.
Note that \begin{equation}
r_{n+1}\leq(r_{n}-p_{\leq n})+(f-p_{n+1})=r_{n}+f-p_{n+1}-p_{\leq n}.\label{eq:ReserveFirefighters}\end{equation}
 This follows since $r_{n}-p_{\leq n}$ is the number of firefighters
placed in $D_{n+2},D_{n+3},\ldots$ for times $1,\ldots,n$, and at
most $f-p_{n+1}$ firefighters are available to be placed in $D_{n+2},D_{n+3},\ldots$
at time $n+1$. Strict inequality occurs if a firefighter is placed
in $D_{k}$ for $k<n+1$ at time $n+1$.

We prove (\ref{eq:HallThmConclusion}) by induction on $n$. For $n=0$,
$\left|B_{0}\right|=1$ holds trivially. We assume the result holds
for $n$, $0\leq n\leq h$, and prove the result for $n+1$. By the
inductive hypothesis, \begin{equation}
\left|B_{n}\right|\geq\begin{cases}
1 & \textrm{if $n=0$,}\\
1+r_{n}+\sum_{i=0}^{n-1}(a_{i}-f) & \textrm{if $1\leq n\leq h$,}\end{cases}\label{eq:InductiveHypothesis1}\end{equation}
 and so by hypotheses~\ref{enu:ConditionKequals0} and \ref{enu:Condition1KH},
\begin{equation}
\left|N(B_{n})\cap D_{n+1}\right|\geq\left|B_{n}\right|+a_{n}.\label{eq:SizeOfBn}\end{equation}
 Thus,\begin{align*}
\left|B_{n+1}\right| & =\left|N(B_{n})\cap D_{n+1}\right|-p_{n+1}-p_{\leq n}\\
 & \geq\left|B_{n}\right|+a_{n}-p_{n+1}-p_{\leq n},\textrm{ by (\ref{eq:SizeOfBn}),}\\
 & \geq1+r_{n}+\sum_{i=0}^{n-1}(a_{i}-f)+a_{n}-p_{n+1}-p_{\leq n},\textrm{ by (\ref{eq:InductiveHypothesis1}),}\\
 & =1+(r_{n}+f-p_{n+1}-p_{\leq n})+\sum_{i=0}^{n-1}(a_{i}-f)+(a_{n}-f)\\
 & \geq1+r_{n+1}+\sum_{i=0}^{n}(a_{i}-f),\textrm{ by (\ref{eq:ReserveFirefighters})}.\end{align*}
 This proves (\ref{eq:HallThmConclusion}) for $0\leq n\leq h+1$.

We now prove (\ref{eq:HallThmConclusion}) for $n\geq h+1$ using
induction on $n$. Note that (\ref{eq:HallThmConclusion}) holds for
$n=h+1$ from above. We thus assume (\ref{eq:HallThmConclusion})
holds for $n\geq h+1$, and we prove the result for $n+1$. By inductive
hypothesis, \begin{equation}
\left|B_{n}\right|\geq1+r_{n}+\sum_{i=0}^{h}(a_{i}-f),\label{eq:InductiveHypothesis2}\end{equation}
 and so by hypothesis \ref{enu:ConditionKgtH}, (\ref{eq:SizeOfBn})
holds for $n>h$. Thus,\begin{align*}
\left|B_{n+1}\right| & =\left|N(B_{n})\cap D_{n+1}\right|-p_{n+1}-p_{\leq n}\\
 & \geq\left|B_{n}\right|+f-p_{n+1}-p_{\leq n},\textrm{ by (\ref{eq:SizeOfBn}),}\\
 & \geq1+r_{n}+\sum_{i=0}^{h}(a_{i}-f)+f-p_{n+1}-p_{\leq n},\textrm{ by (\ref{eq:InductiveHypothesis2}),}\\
 & =1+(r_{n}+f-p_{n+1}-p_{\leq n})+\sum_{i=0}^{h}(a_{i}-f)\\
 & \geq1+r_{n+1}+\sum_{i=0}^{h}(a_{i}-f),\textrm{ by (\ref{eq:ReserveFirefighters}).}\qedhere\end{align*}

This completes the proof of Theorem~\ref{thm:HallTypeCondition}. 
\end{proof}
We now turn our attention to square lattices of dimension three and
higher. It will prove convenient to partition these lattices into
identical subgraphs.

\begin{definition}
The orthants of $\mathbb{R}^{d}$ are the $2^{d}$ regions defined
by the hyperplanes $x_{i}=-1/2$ in $\mathbb{R}^{d}$, $i=1,\ldots,d$.
Let the orthants in $\mathbb{L}^{d}$ be the subsets of vertices that
lie in each orthant of $\mathbb{R}^{d}$. Thus, the $j^{\textrm{th}}$
coordinates of all the vectors in a given orthant of $\mathbb{R}^{d}$
are all non-negative or are all negative, for $j=1,\ldots,d$. Let
$D_{k}^{+}$ denote the vertices of $D_{k}\subseteq\mathbb{L}^{d}$
in the orthant whose elements are all non-negative.

Let $v=(v_{1},v_{2},\ldots,v_{d})$ be an element of $D_{k}\subseteq\mathbb{L}^{d}$.
Let $c_{i}(v)$ denote $v_{i}$, and for a set $A\subseteq D_{k}$
define $A_{r}^{i}=\left\{ v\in A:c_{i}(v)=r\right\} $. Let $v_{\rightarrow i}$
denote $(v_{1},v_{2},\ldots,v_{i}',v_{i+1},\ldots,v_{d})\in D_{k+1}$,
where $v_{i}'=v_{i}+1$ if $v_{i}\geq0$ or $v_{i}'=v_{i}-1$ if $v_{i}<0$.
Thus, $v_{\rightarrow i}$ is in the same orthant as $v$. 
\end{definition}
\begin{lemma}
\label{lem:NeighborsOfDk} In $\mathbb{L}^{d}$ for $d\geq3$, if
$A\subseteq D_{k}$ where $\left|A\right|\geq2d-2$, then $\left|N(A)\cap D_{k+1}\right|\geq\left|A\right|+2d-2$. 
\end{lemma}
\begin{proof}
Given any nonempty set $A\subseteq D_{k}\subseteq\mathbb{L}^{d}$
completely contained in one orthant, we will show that \begin{equation}
\left|N(A)\cap D_{k+1}\right|\geq\left|A\right|+d-1\textrm{, for any $d$}.\label{eq:APlusDMinusOneForAnyAorD}\end{equation}
 We form a set $B\subseteq N(A)\cap D_{k+1}$ in the following way: 
\begin{enumerate}
\item For each $v\in A$, add $v_{\rightarrow1}$ to $B$. 
\item For each $2\leq j\leq d$, let $r_{j}$ be the value of the $j^{\textrm{th}}$
coordinate of elements of $A$ that is greatest in absolute value.
For each $v\in A_{r_{j}}^{j}$, add $v_{\rightarrow j}$ to $B$. 
\end{enumerate}
Each vector added to $B$ in step~1 is unique, and each vector added
to $B$ in step~2 is also unique since the $j^{\textrm{th}}$ coordinate
was chosen to be maximum. Thus, $\left|N(A)\cap D_{k+1}\right|\geq\left|B\right|\geq\left|A\right|+d-1$.

Let $A\subseteq D_{k}\subseteq\mathbb{L}^{d}$. If $A$ is not completely
contained in one orthant, then let $A$ be partitioned as \[
A=A_{1}\cup A_{2}\cup\cdots\cup A_{q},\]
 where each $A_{\ell}$ is in a different orthant $\mathcal{O}_{\ell}$.
By~(\ref{eq:APlusDMinusOneForAnyAorD}), $\left|N(A_{\ell})\cap D_{k+1}\right|\geq\left|A_{\ell}\right|+d-1$.
Note also that the corresponding sets $B_{\ell}$ in the proof above
for $A_{\ell}$ do not overlap since they are in different orthants.
Hence, \begin{align*}
\left|N(A)\cap D_{k+1}\right| & \geq\sum_{\ell=1}^{q}\left|N(A_{\ell})\cap\mathcal{O}_{\ell}\cap D_{k+1}\right|\\
 & \geq\sum_{\ell=1}^{q}\left[\left|A_{\ell}\right|+d-1\right]\\
 & \geq\left|A\right|+2d-2.\end{align*}
 Thus, we may assume that $A$ is completely contained in one orthant,
and, without loss of generality, we assume that all coordinates of
elements of $A$ are non-negative.

We now proceed to prove the lemma by induction on $d$. Let $A\subseteq D_{k}^{+}\subseteq\mathbb{L}^{d}$,
where $\left|A\right|\geq2d-2$. Suppose that $d=3$. Let $n_{i}$
denote the number of nonempty $A_{r}^{i}$, or, equivalently, the
number of distinct $i^{\textrm{th}}$ coordinates of elements of $A$.
Let $i'$ be a coordinate where $n_{i}$ is maximized. We claim that
$n_{i'}\geq3$. If $n_{i'}$ is $1$, then $A$ contains only one
element, which is a contradiction since $\left|A\right|\geq2d-2=4$.
If $n_{i'}$ is $2$, then each coordinate has only two different
values it can assume. However, the sum of the coordinates must remain
$k$. It is straightforward to verify that the maximum number of elements
in $A$ is $3$, which contradicts the fact that $\left|A\right|\geq4$.
Thus, $n_{i'}\geq3$.

For each $r$ where $A_{r}^{i'}$ is nonempty, form a set $\widehat{A_{r}^{i'}}\subseteq D_{k-r}^{d-1}\subseteq\mathbb{L}^{d-1}$
by eliminating the $i'$ coordinate of each element in $A_{r}^{i'}$;
thus, the $\widehat{A_{r}^{i'}}$'s are the parts of $A$ contained
in the slices of $\mathbb{L}^{d}$ taken in direction $i'$. By~(\ref{eq:APlusDMinusOneForAnyAorD}),
$\left|N(\widehat{A_{r}^{i'}})\cap D_{k-r+1}^{d-1}\right|\geq\left|\widehat{A_{r}^{i'}}\right|+d-2$.
For each $v$ in $N(\widehat{A_{r}^{i'}})\cap D_{k-r+1}^{d-1}$, form
an element $\widetilde{v}$ in $N(A_{r}^{i'})\cap D_{k+1}^{d}$ by
inserting $r$ as the $i'$ coordinate. Notice that these elements
are distinct when the $i'$ coordinates are distinct. Let $m$ be
the maximum $r$ such that $A_{r}^{i'}$ is nonempty, or equivalently,
the largest $i'$ coordinate. For each $v\in A_{m}^{i'}$, we also
have $v_{\rightarrow i'}\in N(A)\cap D_{k+1}$, and these vectors
are distinct from any formed above because the $i'$ coordinate is
larger. Thus, \begin{align}
\left|N(A)\cap D_{k+1}\right| & \geq\sum_{r:A_{r}^{i'}\neq\emptyset}\left(\left|A_{r}^{i'}\right|+d-2\right)+\left|A_{m}^{i'}\right|\nonumber \\
 & \geq\left|A\right|+n_{i'}(d-2)+\left|A_{m}^{i'}\right|.\label{eq:BoundWithMaximumNotReplaced}\end{align}
 Since $\left|A_{m}^{i'}\right|\geq1$, (\ref{eq:BoundWithMaximumNotReplaced})
implies that\begin{equation}
\left|N(A)\cap D_{k+1}\right|\geq\left|A\right|+3d-5,\label{eq:ConstructionGivesThreeDMinus5}\end{equation}
 and when $d=3$, \[
\left|N(A)\cap D_{k+1}\right|\geq\left|A\right|+4=\left|A\right|+2d-2.\]

Now suppose that $d>3$. Again let $n_{i}$ denote the number of nonempty
$A_{r}^{i}$, and let $i'$ be a coordinate where $n_{i}$ is maximized.
If $n_{i'}\geq3$, then using the same construction as in the $d=3$
case, we have (\ref{eq:ConstructionGivesThreeDMinus5}), and since
$d>3$, $\left|N(A)\cap D_{k+1}\right|\geq\left|A\right|+2d-2$. If
$n_{i'}=1$, then $A$ contains only one element, which is a contradiction
since $\left|A\right|\geq2d-2\geq4$. We are thus left with the case
$n_{i'}=2$. Let $m$ be the maximum $r$ such that $A_{r}^{i'}$
is nonempty, or equivalently, the largest $i'$ coordinate of elements
of $A$, and let $r'\neq m$ be the minimum value of $r$ where $A_{r}^{i'}$
is nonempty. If $\left|A_{m}^{i'}\right|\geq2$, then using the same
construction as in the $n_{i'}\geq3$ case, we have by (\ref{eq:BoundWithMaximumNotReplaced})
\begin{align*}
\left|N(A)\cap D_{k+1}\right| & \geq\left|A\right|+n_{i'}(d-2)+\left|A_{m}^{i'}\right|\\
 & \geq\left|A\right|+(2d-4)+2,\textrm{ since $\left|A_{m}^{i'}\right|\geq2$,}\\
 & \geq\left|A\right|+2d-2.\end{align*}
 If $\left|A_{m}^{i'}\right|=1$, then we again use the construction
from the $n_{i'}\geq3$ case. However, $\left|\widehat{A_{r'}^{i'}}\right|\geq2d-3$,
so by induction, $\left|N(\widehat{A_{r'}^{i'}})\cap D_{k-r'+1}^{d-1}\right|\geq\left|\widehat{A_{r'}^{i'}}\right|+2d-4$.
Here, the notation $D_{z}^{d-1}$ means the set $D_{z}\subseteq\mathbb{L}^{d-1}$,
emphasizing the dimension of $\mathbb{L}^{d-1}$. For each $v$ in
$N(\widehat{A_{r'}^{i'}})\cap D_{k-r'+1}^{d-1}$, form an element
$\widetilde{v}$ in $N(A_{r'}^{i'})\cap D_{k+1}^{d}$ by inserting
$r'$ as the $i'$ coordinate. Additionally, for the single vector
$v\in A_{m}^{i'}$ and $1\leq j\leq d$, $v_{\rightarrow j}\in N(A)\cap D_{k+1}$,
and these vectors are distinct from those formed above because the
$i'$ coordinate is larger. Thus, \begin{align*}
\left|N(A)\cap D_{k+1}\right| & \geq\left(\left|A_{r'}^{i'}\right|+2d-4\right)+d\\
 & =\left|A\right|+3d-3,\textrm{ since $\left|A_{r'}^{i'}\right|=\left|A\right|+1$,}\\
 & \geq\left|A\right|+2d-2,\textrm{ since $d>3$.}\qedhere\end{align*}

\end{proof}
Lemma~\ref{lem:NeighborsOfDk} provides the long-term growth of the
front, Condition 3 needed for Theorem~\ref{thm:HallTypeCondition}.
The next lemma gives the complementary requirements.

\begin{lemma}
\label{lem:NeighborsOfD1} In $\mathbb{L}^{d}$ for $d\geq3$, if
$A\subseteq D_{1}$ where $\left|A\right|\geq2$, then $\left|N(A)\cap D_{2}\right|\geq\left|A\right|+4d-6$. 
\end{lemma}
\begin{proof}
Let $A\subseteq D_{1}\subseteq\mathbb{L}^{d}$ where $\left|A\right|\geq2$.
Every vector $v\in A$ is of the form $(0,\ldots,x_{i},\ldots,0)$,
where $x_{i}=\pm1$. Each vector $v$ in $A$ has $2(d-1)$ neighbors
in $D_{2}$ formed by replacing each of the zero coordinates in $v$
with $\pm1$, and one neighbor formed by replacing $1$ in the $i^{\textrm{th}}$
coordinate with $2$ or replacing $-1$ with $-2$. If $v$ and $v'$
are vectors of $A$ with nonzero entries in different coordinates,
then $v$ and $v'$ share exactly one neighbor in $D_{2}$. If $v$
and $v'$ have nonzero entries in the same coordinate, then $v$ and
$v'$ share no neighbors in $D_{2}$. Thus, \begin{align*}
\left|N(A)\cap D_{2}\right| & \geq\left|A\right|\left(2(d-1)+1\right)-{{\left|A\right| \choose 2}}\\
 & =\left|A\right|\left(2d-\frac{\left|A\right|}{2}-\frac{1}{2}\right)\\
 & \geq\left|A\right|+\left|A\right|\left(2d-\frac{\left|A\right|}{2}-\frac{3}{2}\right).\end{align*}
 It is straightforward to use calculus to verify that \[
\left|A\right|\left(2d-\frac{\left|A\right|}{2}-\frac{3}{2}\right)\geq4d-6,\]
 where $d\geq3$ and $2\leq\left|A\right|\leq2d$, and so \[
\left|N(A)\cap D_{2}\right|\geq\left|A\right|+4d-6.\qedhere\]

\end{proof}
\begin{theorem}
In $\mathbb{L}^{d}$, $2d-1$ firefighters are needed to contain an
outbreak of fire starting at a single vertex. 
\end{theorem}
\begin{proof}
Since $\mathbb{L}^{d}$ is vertex transitive, we may assume that the
root vertex where the fire outbreak starts is the origin. We use Theorem~\ref{thm:HallTypeCondition}
with $f=2d-2$, $h=1$, $a_{0}=2d-1$, and $a_{1}=4d-6$. The one
element set $D_{0}$ has $2d$ neighbors in $D_{1}$ so hypothesis~\ref{enu:ConditionKequals0}
of Theorem~\ref{thm:HallTypeCondition} holds, Lemma~\ref{lem:NeighborsOfD1}
shows hypothesis~\ref{enu:Condition1KH} of Theorem~\ref{thm:HallTypeCondition}
holds for $k=1$, and Lemma~\ref{lem:NeighborsOfDk} shows hypothesis~\ref{enu:ConditionKgtH}
holds for $k>1$. By Theorem~\ref{thm:HallTypeCondition}, $2d-2$
firefighters are insufficient to contain an outbreak starting at the
origin. 
\end{proof}
Fogarty also showed in \cite{Fogarty} that two firefighters suffice
in $\mathbb{L}^{2}$ to contain any finite outbreak of fire where
an arbitrarily large but finite number of vertices are initially on
fire. However, we prove for $\mathbb{L}^{d}$ where $d\geq3$ that
for any fixed number $f$ of firefighters, there is a finite outbreak
of fire in which $f$ firefighters per time step are insufficient
to contain the outbreak.

First we establish the following lemma. Essentially, the lemma says
that if we have a {}``front'' of $x$ elements, then it will grow
outwards by at least $\Omega(\sqrt{x})$ in the next time step.

\begin{lemma}
\label{lem:GrowthOfDkinL3}Let $f$ be any positive integer. If $A\subseteq D_{k}^{+}\subseteq\mathbb{L}^{3}$
where $\left|A\right|>\frac{1}{2}(f-1)(f-2)$, then $\left|N(A)\cap D_{k+1}^{+}\right|\geq\left|A\right|+f$. 
\end{lemma}
\begin{proof}
Let $A\subseteq D_{k}^{+}\subseteq\mathbb{L}^{3}$ be a set where
$\left|A\right|>\frac{1}{2}(f-1)(f-2)$. The elements of $B:=\left\{ v_{\rightarrow1}:v\in A\right\} $
are distinct vertices in $N(A)\cap D_{k+1}^{+}$, and the set $B$
has cardinality equal to $|A|$. Therefore, it suffices to show that
if $\left|A\right|>\frac{1}{2}(f-1)(f-2)$, then there are at least
$f$ distinct elements of the form $v_{\rightarrow j}$ which are
not elements of $B$, where $v\in A$ and $j\in\{2,3\}$.

Let $m$ be the largest first coordinate of elements of $A$, and
let $t$ be the smallest first coordinate of elements of $A$. Recall
that the sets $A_{r}^{1}$, $r=t,t+1,\ldots,m$, partition $A$. Let
$\sigma_{r}$ equal $\left|A_{r}^{1}\right|$, so that $\sum_{r=t}^{m}\sigma_{r}=|A|$.
Note that $\sigma_{t},\sigma_{m}>0$.

Suppose some $\sigma_{r}$ is equal to zero, where $t<r<m$. Then
$A$ is partitioned into the sets $A_{1}$ consisting of all elements
of $A$ with first coordinate greater than $r$ and $A_{2}$ consisting
of all elements of $A$ with first coordinate less than $r$. Clearly
$N(A_{1})\cap N(A_{2})\cap D_{k+1}^{+}=\varnothing$. Define $A_{1}^{\prime}:=\left\{ v_{\rightarrow2}:v\in A_{1}\right\} $
and $A_{2}^{\prime}:=\left\{ v_{\rightarrow1}:v\in A_{2}\right\} $,
so that $A_{1}^{\prime}$ and $A_{2}^{\prime}$ are subsets of $D_{k+1}^{+}$.
Since $A_{1}'$ is simply a translate of $A_{1}$ by $1$ in the first
coordinate, $N(A_{1}')\cap D_{k+2}^{+}$ is a translate of $N(A_{1})\cap D_{k+1}^{+}$
by $1$ in the first coordinate. Similarly, $N(A_{2}')\cap D_{k+2}^{+}$
is a translate of $N(A_{2})\cap D_{k+1}^{+}$ by $1$ in the second
coordinate. Thus, we have that \begin{align*}
\left|N(A_{1}^{\prime}\cup A_{2}^{\prime})\cap D_{k+2}^{+}\right| & \leq\left|N(A_{1}')\cap D_{k+2}^{+}\right|+\left|N(A_{2}')\cap D_{k+2}^{+}\right|\\
 & =\left|N(A_{1})\cap D_{k+1}^{+}\right|+\left|N(A_{2})\cap D_{k+1}^{+}\right|\\
 & =\left|N(A)\cap D_{k+1}^{+}\right|,\end{align*}
 where the last equality follows since $N(A_{1})\cap D_{k+1}^{+}$
and $N(A_{2})\cap D_{k+1}^{+}$ do not intersect. However, $A_{1}^{\prime}\cup A_{2}^{\prime}$
has the same size as $A$, but the separation between the largest
first coordinate of elements of $A_{1}^{\prime}\cup A_{2}^{\prime}$
and the smallest first coordinate of $A_{1}^{\prime}\cup A_{2}^{\prime}$
is less than $m-t$. Therefore, by induction on $m-t$ we reduce to
the case where no $\sigma_{r}$ is equal to zero, \emph{i.e.}, there
is an element of $A$ with first coordinate $r$ for every $t\leq r\leq m$.

Consider the sets $S_{r}=\left\{ v_{\rightarrow j}:v\in A_{r}^{1},j\in\{2,3\}\right\} \subseteq N(A)\cap D_{k+1}^{+}$.
Observe that the cardinality of $S_{r}$ is at least $\sigma_{r}+1$.
Clearly all $S_{r}$ are disjoint, since all elements of $S_{r}$
have first coordinate $r$. The elements of $S_{t}$ have $t$ as
their first coordinate, while all elements of $B$ have first coordinates
at least $t+1$, so no elements of $S_{t}$ are in $B$. Furthermore,
for all $r>t$, if an element of $S_{r}$ is in $B$, then by considering
its first coordinate, the element must be in the set $\left\{ v_{\rightarrow1}:v\in A_{r-1}^{1}\right\} $.
In particular, this set has size $\sigma_{r-1}$. If $\sigma_{r}+1>\sigma_{r-1}$,
then there are at least $\sigma_{r}+1-\sigma_{r-1}$ elements in $S_{r}$
not in $B$. Therefore, the number of elements in $N(A)\cap D_{k+1}^{+}$
that are not in $B$ is bounded below by \begin{equation}
g(\sigma):=\sum_{r=t}^{m}\text{max}\:(0,\sigma_{r}+1-\sigma_{r-1}),\label{eq:sigma-excess}\end{equation}
 with the convention that $\sigma_{t-1}=0$.

Now take any nonzero sequence $\sigma_{t},\sigma_{t+1},\ldots,\sigma_{m}$.
We claim that if $g(\sigma)<f$, then $\sum_{r=t}^{m}\sigma_{r}\leq\frac{1}{2}(f-1)(f-2)$,
which would complete the proof of the theorem. Suppose we have some
sequence $\sigma_{t},\sigma_{t+1},\ldots,\sigma_{m}$ with $g(\sigma)<f$.
First, suppose that there exists some $r>t$ where $\sigma_{r}\geq\sigma_{r-1}$.
Then adding $1$ to $\sigma_{r-1}$ decreases the $r$-th term of
(\ref{eq:sigma-excess}) by $1$, possibly adds $1$ to the $(r-1)$-st
term, and leaves all other terms unchanged; in particular, it does
not increase the value of $g(\sigma)$ and increases $\sum\sigma_{r}$.
Therefore, we can reduce to the case where $\sigma$ is strictly decreasing.

Next, suppose we have $\sigma_{r}<\sigma_{r-1}-1$ for some $t<r\leq m$.
Then adding $1$ to $\sigma_{r}$ leaves all terms of (\ref{eq:sigma-excess})
unchanged. Similar to before, this operation does not change $g(\sigma)$,
while increasing $\sum\sigma_{r}$. Doing this repeatedly, we reduce
to the case where \begin{equation}
\sigma_{r-1}=\sigma_{r}+1\label{eq:FinalSigmaSeq}\end{equation}
 for all $t<r\leq m$. However, this case is easy to evaluate; each
term in (\ref{eq:sigma-excess}) is zero except the $r=t$ term, which
is equal to $\sigma_{t}+1$. Since $g(\sigma)=\sigma_{t}+1<f$, $\sigma_{t}<f-1$.
Since $\sigma_{m}>0$, $\sum_{r=t}^{m}\sigma_{r}$ is at most the
sum of the first $f-2$ positive integers. Thus, \[
\sum_{r=t}^{m}\sigma_{r}\leq\frac{1}{2}(f-1)(f-2).\qedhere\]

\end{proof}
This allows us to prove the following theorem.

\begin{theorem}
\label{thm:CanNotContainAllFiniteOutbreaks}For any dimension $d\geq3$
and any fixed positive integer $f$, $f$ firefighters per time step
are not sufficient to contain all finite outbreaks in $\mathbb{L}^{d}$. 
\end{theorem}
\begin{proof}
Since $\mathbb{L}^{3}$ is contained in $\mathbb{L}^{d}$ for $d\geq3$,
it suffices to prove the statement for $d=3$. We consider an initial
outbreak consisting of all of $D_{k}^{+}$ for $k$ large enough so
that $\left|D_{k}^{+}\right|>\frac{1}{2}(f-1)(f-2)$. To show that
$f$ firefighters are insufficient to contain this outbreak, we will
construct a related graph that captures the essential disease dynamics
and then invoke Theorem~\ref{thm:HallTypeCondition}. Let $G$ be
the subgraph of $\mathbb{L}^{3}$ induced by vertices with non-negative
coordinates that are distance at least $k$ from the origin. Let $G'$
be the graph formed from $G$ by identifying all of the vertices in
$D_{k}^{+}$ as a single vertex $r$. An edge exists between vertices
$x$ and $y$ in $G'$ if $xy$ is an edge in $G$ or if $x=r$ and
$y\in N_{G}(D_{k}^{+})$. Let $D_{i}'$ denote the set of vertices
in $G'$ that are distance $i$ from the root $r$. By Lemma~\ref{lem:GrowthOfDkinL3},
\[
\left|N(D_{k}^{+})\cap D_{k+1}^{+}\right|\geq\left|D_{k}^{+}\right|+f>\frac{1}{2}(f-1)(f-2)+f,\]
 and so \[
\left|N(r)\cap D_{1}'\right|>\left(\left|D_{0}'\right|-1\right)+\frac{1}{2}(f-1)(f-2)+f.\]
 If $A'\subseteq D_{i}'$, where $i>0$ and $\left|A'\right|>\frac{1}{2}(f-1)(f-2)$,
then $A'$ corresponds to a set $A\subseteq D_{k+i}^{+}$ and by Lemma~\ref{lem:GrowthOfDkinL3},
\[
\left|N(A)\cap D_{k+i+1}^{+}\right|\geq\left|A\right|+f,\]
 and hence \[
\left|N(A')\cap D_{i+1}'\right|\geq\left|A'\right|+f.\]
 By Theorem~\ref{thm:HallTypeCondition} with $h=0$, and $a_{0}=\frac{1}{2}(f-1)(f-2)+f$,
$f$ firefighters are insufficient to contain an outbreak starting
at $r$ in $G'$, and hence $f$ firefighters are insufficient to
contain an outbreak consisting of all of $D_{k}^{+}$ in $\mathbb{L}^{3}$. 
\end{proof}
The essential problem here is that for $d\geq3$, the boundary of
an outbreak grows faster than the constant number of firefighters
deployed at a given time step. Indeed, in dimension $d$, the boundary
grows as a polynomial of degree $d-2$. This motivates the following
ambitious conjecture.

\begin{conjecture}
Suppose that $f(t)$ is a function on $\mathbb{N}$ with the property
that $\frac{f(t)}{t^{d-2}}$ goes to $0$ as $t$ gets large. Then
there exists some outbreak on $\mathbb{L}^{d}$ which cannot be contained
by deploying $f(t)$ firefighters at time $t$. 
\end{conjecture}
A weaker conjecture would require $f(t)$ to be a polynomial.

Lemma~\ref{lem:GrowthOfDkinL3} also allows us to resolve another
conjecture of Wang and Moeller in \cite{WangMoeller}. They conjectured
that as $n$ gets large, the proportion of elements in the three-dimensional
grid $P_{n}\times P_{n}\times P_{n}$ which can be saved by using
one firefighter per time step when an outbreak at one vertex occurs
goes to $0$. We prove this conjecture in the following

\begin{theorem}
Let $v$ be any vertex of $P_{n}\times P_{n}\times P_{n}$, for $n\geq1$.
Then the maximum number of vertices which can be saved by deploying
one firefighter per time step with an initial outbreak at $v$ grows
at most as $O(n^{2})$. In particular, the proportion of vertices
which can be saved goes to $0$ as $n$ gets large. 
\end{theorem}
\begin{proof}
We prove the theorem in the case $v=(0,0,0)$. The general statement
easily follows by splitting $P_{n}\times P_{n}\times P_{n}$ into
orthants with apex $v$. We actually prove a stronger statement. Consider
the graph $G$ induced from the lattice $\mathbb{L}^{3}$ by vertices
with non-negative coordinates and distance at most $3n$ from the
origin $v$. We prove the theorem for the graph $G$. Note that $G$
contains $P_{n}\times P_{n}\times P_{n}$ as an induced subgraph.

We claim that $\left|B_{t}\right|-r_{t}\geq\frac{t^{2}+t+2}{2}$ for
all $t$ regardless of what firefighter placements are made. Since
there are ${{t+2 \choose 2}}=\frac{t^{2}+3t+2}{2}$ vertices in $D_{t}^{+}$,
this statement is saying that at the end of the $t^{\textrm{th}}$
time step the number of reserve firefighters together with the unburned
vertices (including defended vertices) in $D_{t}^{+}$ cannot exceed
$t$. By considering time up to $t=3n$, when all vertices have had
a chance to be burned, at most $1+2+\ldots+3n=O(n^{2})$ vertices
are unburned. This implies the same conclusion for $P_{n}\times P_{n}\times P_{n}$.

The proof of the claim is by induction. At the end of the $0^{\textrm{th}}$
time step, there are no reserve firefighters, and one vertex in $D_{1}$
is burned; the difference is $1-0=1\geq1=\frac{0^{2}+0+2}{2}$ as
desired.

Suppose $t\geq0$, and suppose that the statement is true for $t$.
Then \begin{equation}
\left|B_{t}\right|\geq\frac{t^{2}+t+2}{2}>\frac{1}{2}t(t+1).\label{eq:InductiveHypothesisForPn3}\end{equation}
 Let $f=t+2$. By Lemma~\ref{lem:GrowthOfDkinL3}, \begin{equation}
\left|N(B_{t})\cap D_{t+1}^{+}\right|\geq\left|B_{t}\right|+f.\label{eq:SizeofBnp1ForPn3}\end{equation}
 As in the proof of Theorem~\ref{thm:HallTypeCondition}, let $p_{t+1}$
denote the number of firefighters placed in $D_{t+1}^{+}$ at time
$t+1$, and let $p_{\leq t}$ denote the number of reserve firefighters
placed in $D_{t+1}^{+}$ during time steps $1,\ldots,t$. Thus,\begin{align*}
\left|B_{t+1}\right|-r_{t+1} & =\left[\left|N(B_{t})\cap D_{t+1}^{+}\right|-p_{t+1}-p_{\leq t}\right]-r_{t+1}\\
 & \geq\left|N(B_{t})\cap D_{t+1}^{+}\right|-r_{t}-1,\textrm{ by (\ref{eq:ReserveFirefighters}),}\\
 & \geq\left|B_{t}\right|+f-r_{t}-1,\textrm{ by (\ref{eq:SizeofBnp1ForPn3}),}\\
 & \geq\frac{t^{2}+t+2}{2}+(t+2)-1,\textrm{ by (\ref{eq:InductiveHypothesisForPn3}),}\\
 & \geq\frac{(t+1)^{2}+(t+1)+2}{2}.\end{align*}
 Hence the claim follows. 
\end{proof}
In practice, one can ensure when an outbreak starts at $(0,0,0)$
that $t$ vertices in $D_{t}^{+}$ are unburned at time $t$. However,
because the fire doubles back on itself, it is unclear that one can
actually save a quadratic number of vertices. Wang and Moeller exhibit
the construction of building a {}``fire wall'' by defending all
of the vertices at distance $k$ from $(n,n,n)$. In order for this
to be effective, we must be able to cover all $\frac{(k+1)(k+2)}{2}$
such vertices in the $3n-k$ time steps it takes the fire to reach
this hyperplane. This yields $k=O(\sqrt{n})$. The number of vertices
saved is the number of vertices at distance $k$ or less from $(n,n,n)$,
which is $\frac{(k+1)(k+2)(k+3)}{6}$. This is $O(k^{3})=O(n^{3/2})$.
Therefore, the optimal number of vertices saved given an initial outbreak
at $(0,0,0)$ in the grid graph $P_{n}\times P_{n}\times P_{n}$ when
deploying one firefighter per time step is between $O(n^{3/2})$ and
$O(n^{2})$.

\section{\label{sec:TwoDgrid}Two Dimensional Square Grid}

According to Wang and Moeller in \cite{WangMoeller}, Hartnell, Finbow,
and Schmeisser first proved that an outbreak of fire in $\mathbb{L}^{2}$
starting at a single vertex can be contained using two firefighters
per time step. Their sequence of firefighter placements contained
the outbreak at the end of $11$ time steps. Wang and Moeller showed
that the disease cannot be contained at the end of $7$ time steps
when using two firefighters per time step and presented a sequence
of firefighter placements that attains this minimum. Their sequence
allows $18$ vertices to be burned. Surprisingly, Wang and Moeller
do not comment on whether their solution attains the minimum number
of burned vertices. In fact, $18$ is the minimum number of burned
vertices, and we prove this using integer programming. The same technique
also gives a computer proof of Wang and Moeller's result that at least
$8$ time steps are needed. Their proof relies heavily on case analysis.

The tightness in the following theorem is due to Wang and Moeller~\cite{WangMoeller}.

\begin{theorem}
\label{thm:EighteenVerticesBurnedinL2}In $\mathbb{L}^{2}$, if an
outbreak of fire starts an a single vertex, then when using two firefighters
per time step at least $18$ vertices are burned. This bound is tight. 
\end{theorem}
\begin{proof}
We formulate an integer program using the boolean variables $b_{x,t}$
and $d_{x,t}$. The variable $b_{x,t}$ is $1$ if and only if vertex
$x$ is burned at or before time $t$, and $d_{x,t}$ is $1$ if and
only if $x$ is defended at or before time $t$. We wish to minimize
the total number of vertices that become burned. For the integer program
to be implementable with a finite number of variables and constraints,
we restrict the graph to $L=\{(x,y)\in\mathbb{L}^{2}:\textrm{$\left|x\right|\leq\ell$ and $\left|y\right|\leq\ell$}\}$
and $0\leq t\leq T$, where $\ell$ and $T$ are chosen to be sufficiently
large that the fire never reaches the boundary and is completely contained
by time $T$. In the actual computations performed, $\ell=6$ and
$T=9$ proved sufficient. We choose $T>8$ to ensure that the fire
is actually contained and does not grow in the last time step.

The integer program is \begin{align}
\textrm{minimize} & \sum_{x\in L}b_{x,T}\nonumber \\
\textrm{subject to:} & b_{x,t}+d_{x,t}-b_{y,t-1}\geq0,\textrm{ for all $x\in L$, $y\in N(x)$, and $1\leq t\leq T$,}\label{eq:Grid:FireSpreads}\\
 & b_{x,t}+d_{x,t}\leq1,\textrm{ for all $x\in L$ and $1\leq t\leq T$,}\label{eq:Grid:DontDefendBurntVertex}\\
 & b_{x,t}-b_{x,t-1}\geq0,\textrm{ for all $x\in L$ and $1\leq t\leq T$,}\label{eq:Grid:StayBurnt}\\
 & d_{x,t}-d_{x,t-1}\geq0,\textrm{ for all $x\in L$ and $1\leq t\leq T$,}\label{eq:Grid:StayDefended}\\
 & \sum_{x\in L}\left(d_{x,t}-d_{x,t-1}\right)\leq2,\textrm{ for $1\leq t\leq T$,}\label{eq:Grid:TwoFirefightersPerTimeStep}\\
 & b_{x,0}=\begin{cases}
1 & \textrm{if $x$ is the origin,}\\
0 & \textrm{otherwise,}\end{cases}\textrm{ for all $x\in L$,}\label{eq:Grid:BurnInitialCondition}\\
 & d_{x,0}=0,\textrm{ for all $x\in L$,}\label{eq:Grid:DefendedInitialCondition}\\
 & b_{x,t},d_{x,t}\in\{0,1\},\textrm{ for all $x\in L$ and $0\leq t\leq T$.}\label{eq:Grid:BinaryConstraint}\end{align}
 Condition~(\ref{eq:Grid:FireSpreads}) enforces the spread of the
fire while respecting vertices defended by a firefighter. Note that
vertices can spontaneously combust, catching fire, but the minimization
of the objective function ensures that this does not happen in the
optimal solution. Condition~(\ref{eq:Grid:DontDefendBurntVertex})
prevents a firefighter from defending a burnt vertex, while conditions
(\ref{eq:Grid:StayBurnt}) and (\ref{eq:Grid:StayDefended}) ensure
that once a vertex is burnt or defended, it stays in that state. Condition~(\ref{eq:Grid:TwoFirefightersPerTimeStep})
only allows two firefighters per time step. Conditions (\ref{eq:Grid:BurnInitialCondition})
and (\ref{eq:Grid:DefendedInitialCondition}) give the initial conditions
at time $t=0$, and condition~(\ref{eq:Grid:BinaryConstraint}) makes
the program a binary integer program.

The integer program was solved in about $1.83$ hours using the GNU
Linear Programming Kit~\cite{GLPK} running on a Pentium IV 2.6GHz
processor, and $18$ was the minimum number of burnt vertices at time
$t=9$. %
\begin{figure}
\begin{center}\scalebox{0.9}{\begin{picture}(0,0)\includegraphics{rectangulargrid.eps}\end{picture}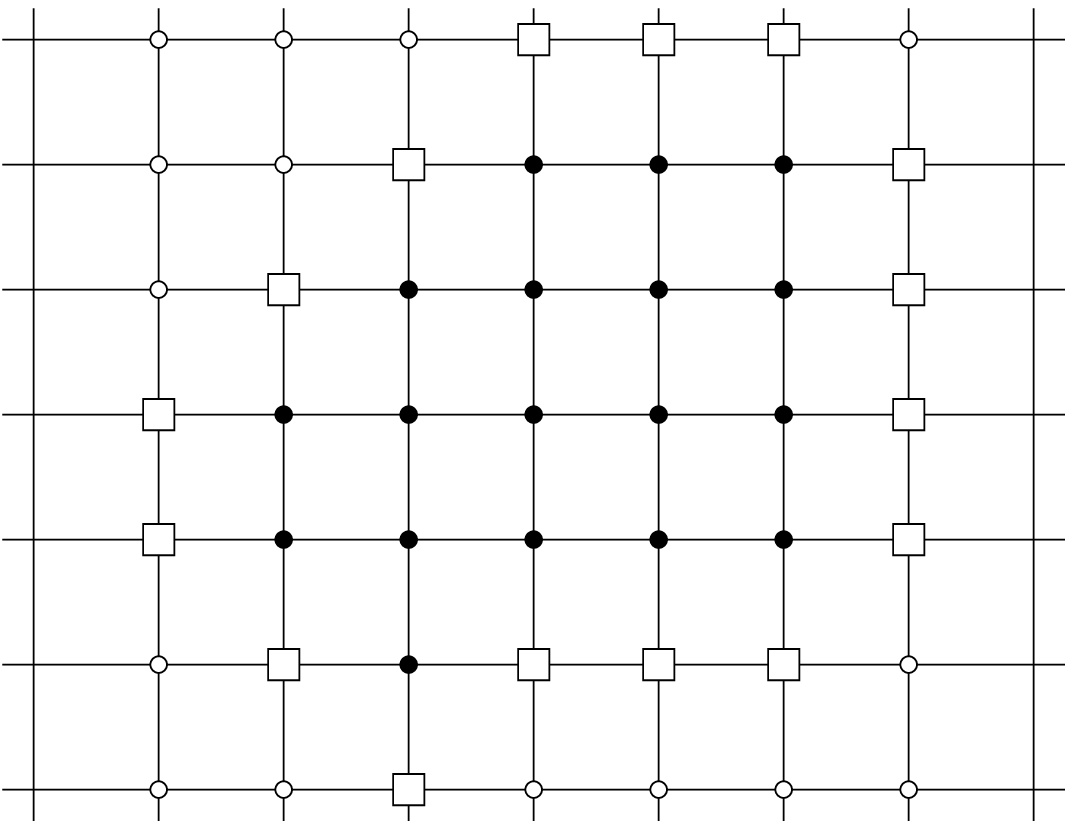}\end{center}

\caption{\label{cap:ContainmentInL2} Optimal solution of the integer program
used in the proof of Theorem~\ref{thm:EighteenVerticesBurnedinL2}.
The fire outbreak starts at time $0$ at the root, and then spreads
to the black vertices at the times written next to the vertices. The
square firefighters $a_{i}$ are placed at time $i$. This placement
of two firefighters per time step in $\mathbb{L}^{2}$ completely
contains the outbreak in $8$ time steps, allowing only the minimum
number of $18$ burned vertices.}
\end{figure}
 Figure~\ref{cap:ContainmentInL2} shows the minimum solution. The
fire was completely contained and thus did not reach the sides of
$L$. Also note that the solution presented by Wang and Moeller in
\cite{WangMoeller} also allows only $18$ burnt vertices but is slightly
different from the solution presented here. 
\end{proof}
\begin{lemma}
\label{lem:NoCoordinateReaches5}If an outbreak of fire in $\mathbb{L}^{2}$
is contained by $14$ defended vertices and $(x,y)$ is a burnt vertex,
then $\left|x\right|\leq5$ and $\left|y\right|\leq5$. 
\end{lemma}
\begin{proof}
Suppose that $(x,y)$ is a burnt vertex, and, without loss of generality,
that $x>5$. Since $(x,y)$ is burnt, there is a path $v_{0}=(x,y),v_{1},v_{2},\ldots,v_{t}=(0,0)$
from $(x,y)$ to the origin consisting of burnt vertices. For each
$0\leq a\leq6$, there is a vertex $v_{\rho(a)}$ such that the first
coordinate of $v_{\rho(a)}$ is $a$. Since the fire is contained,
there must be a defended vertex above and below each of these seven
vertices, and there must be at least one defended vertex with first
coordinate less than $0$ and one with first coordinate greater than
$x$. But this requires $16$ defended vertices, resulting in a contradiction. 
\end{proof}
\begin{theorem}
[Wang and Moeller]In $\mathbb{L}^{2}$, if an outbreak of fire starts
at a single vertex, then the fire cannot be contained at the end of
$7$ time steps when using two firefighters per time step. Thus, at
least $8$ time steps are needed to contain the fire, and this bound
is tight. 
\end{theorem}
\begin{proof}
We use a similar integer program to the one used in the proof of Theorem~\ref{thm:EighteenVerticesBurnedinL2}.
By Lemma~\ref{lem:NoCoordinateReaches5}, if the outbreak can be
contained after $7$ time steps, then no burnt vertex will have either
coordinate equaling $6$ in absolute value. We thus use the finite
grid $L$ where $\ell=6$, and we use the objective function \[
\textrm{minimize }\sum_{\substack{x=(a,b)\in L\\
\left|a\right|=6\textrm{ or }\left|b\right|=6}
}b_{x,T}.\]
 If the disease can be contained after $7$ time steps, then the optimal
value of the objective function will be $0$. All of the conditions
from the previous integer program are included except condition~(\ref{eq:Grid:TwoFirefightersPerTimeStep})
is changed to \begin{equation}
\sum_{x\in L}\left(d_{x,t}-d_{x,t-1}\right)\leq\begin{cases}
2 & \textrm{for $1\leq t\leq7$,}\\
0 & \textrm{for $8\leq t\leq T$.}\end{cases}\label{eq:Grid:NoExtraTimeSteps}\end{equation}
 This prevents firefighters from being used after $7$ time steps.

The integer program with $T=9$ was solved in about $40$ minutes
using the GNU Linear Programming Kit running on a Pentium M 900MHz
processor. The minimum value was $1$, meaning that in every feasible
solution, the fire burned a vertex with one coordinate equaling $6$
in absolute value. This contradicts Lemma~\ref{lem:NoCoordinateReaches5},
and so at least $8$ time steps are needed to contain an outbreak
in $\mathbb{L}^{2}$ when using two firefighters per time step. 
\end{proof}

\section{Future Work}

There are many avenues for future work in models of responses to fire
and disease spread. For infinite graphs, we can ask the same question
as for the infinite square grids: What is the minimum number of firefighters
needed per time step so that only a finite number of vertices are
burned? Percolation is a related topic whose methods may also apply
here.

From the viewpoint of an arsonist or bioterrorist, one would like
to find the most vulnerable vertices in a graph $G$. A vertex $v$
is most vulnerable if a fire outbreak starting at $v$ burns the most
vertices of $G$ given an optimal firefighter response. Can the most
vulnerable vertices in a graph be determined without knowing the optimal
firefighter response? Perhaps they could then be preemptively defended.
From the viewpoint of a network architect, we would like to design
graphs that are resistant to such attacks. Similar questions can also
be asked if there are $k$ initial outbreaks of fire.

Finally, MacGillivray and Wang~\cite{MacGillivrayWang} observed
that the firefighter problem can be viewed as a one-player game. Suppose
that the fire has a choice, too: the fire can only spread to $d$
neighbors each time step. This forms a two-player game. What strategy
should the firefighters use to minimize the number of burned vertices?

\section*{Acknowledgements}

The second author thanks Fred Roberts and James Abello for discussions
and encouragement and Mike Dinitz for providing references.


\begin{thebibliography}{1}
\bibitem{BurkeEpstein}J.~M.~Epstein, D.~A.~T.~Cummings, S.~Chakravarty, R.~M.~Singa,
and D.~ S.~Burke, Toward a Containment Strategy for Smallpox Bioterror:
An Individual-Based Computational Approach, The Brookings Institute
Center on Social and Economic Dynamics Working Paper No. 31, December
2002. 
\bibitem{EpiSIMS}S.~Eubank, H.~Guclu, V.~S.~A.~Kumar, M.~V.~Marathe, A.~Srinivasan,
Z.~Toroczkai, N.~Wang, and the EpiSims Team, \texttt{http://episims.lanl.gov},
April 22, 2004. 
\bibitem{Fogarty}P. Fogarty, \emph{Catching the Fire on Grids}, M.Sc. Thesis, Department
of Mathematics, University of Vermont, 2003. 
\bibitem{GLPK}GNU Linear Programming Kit, \texttt{http://www.gnu.org/software/glpk/glpk.html}. 
\bibitem{Hartnell}B.~Hartnell, Firefighter! An Application of Domination, presentation,
Twentieth Conference on Numerical Mathematics and Computing, University
of Manitoba in Winnipeg, Canada, Sept. 1995. 
\bibitem{MacGillivrayWang}G. MacGillivray and P. Wang, On the Firefighter Problem, \emph{J.
Combin. Math. Combin. Comput.}, 47 (2003), 83-96. 
\bibitem{ReadKeeling}J. M. Read and M. J. Keeling, Disease Evolution on Networks: The Role
of Contact Structure, \emph{Proc. Roy. Soc. Lond. B}, 270 (2003),
699-708. 
\bibitem{WangMoeller}P. Wang and S. A. Moeller, Fire Control on Graphs, \emph{J. Combin.
Math. Combin. Comput.}, 41 (2002), 19-34.\end{thebibliography}
\end{document}